\documentclass[number,citesort,number,citesort,seceqn,dvips]{arxbj}
\usepackage{graphicx}
\usepackage{breakurl}

\aid{0}
\volume{17}
\issue{2}
\pubyear{2011}
\firstpage{592}
\lastpage{608}
\doi{10.3150/10-BEJ295}

\begin{document}

\newcommand{\Bor}{\operatorname{Bor}}
\newcommand{\card}{\operatorname{card}}
\newcommand{\CI}{\operatorname{CI}}
\newcommand{\conv}{\operatorname{conv}}
\newcommand{\Cov}{\operatorname{Cov}}
\newcommand{\CS}{\operatorname{CS}}
\newcommand{\cx}{\operatorname{cx}}
\newcommand{\ess}{\operatorname{ess}}
\newcommand{\Ex}{\operatorname{Ex}}
\newcommand{\hr}{\operatorname{hr}}
\newcommand{\icx}{\operatorname{icx}}
\newcommand{\IE}{\operatorname{IE}}
\newcommand{\iintegr}{\operatorname{iint}}
\newcommand{\imax}{\operatorname{imax}}
\newcommand{\integr}{\operatorname{I}}
\newcommand{\Leb}{\operatorname{Leb}}
\newcommand{\lcm}{\operatorname{lcm}}
\newcommand{\Lo}{\operatorname{Lo}}
\newcommand{\Pa}{\operatorname{Pa}}
\newcommand{\per}{\mathbf{per}}
\newcommand{\priv}{\operatorname{priv}}
\newcommand{\publ}{\operatorname{publ}}
\newcommand{\refin}{\mathrm{ref.}}
\newcommand{\SUPREMUM}{\operatorname{S}}
\newcommand{\schur}{\operatorname{sc}}
\newcommand{\st}{\operatorname{st}}
\newcommand{\unif}{\operatorname{unif}}
\newcommand{\val}{\operatorname{val}}
\newcommand{\Var}{\operatorname{Var}}
\newcommand{\diff}{\mathrm{d}}
\newcommand{\eqref}[1]{(\ref{#1})}

\newtheorem{theorem}{Theorem}[section]
\newtheorem{corollary}[theorem]{Corollary}
\newtheorem{lemma}[theorem]{Lemma}
\newtheorem{lemmasss}{Lemma}
\newtheorem{proposition}[theorem]{Proposition}
\newtheorem{acknowledgement}[theorem]{Acknowledgement}
\newtheorem{algorithm}[theorem]{Algorithm}
\newtheorem{axiom}[theorem]{Axiom}
\newtheorem{case}[theorem]{Case}
\newtheorem{claim}[theorem]{Claim}
\newtheorem{fact}[theorem]{Fact}
\newtheorem{conclusion}[theorem]{Conclusion}
\newtheorem{condition}[theorem]{Condition}
\newtheorem{conjecture}[theorem]{Conjecture}
\newtheorem{criterion}[theorem]{Criterion}
\newtheorem{definition}[theorem]{Definition}
\newremark{example}[theorem]{Example}
\newtheorem{exercise}[theorem]{Exercise}
\newtheorem{notation}[theorem]{Notation}
\newtheorem{problem}[theorem]{Problem}
\newtheorem{remark}[theorem]{Remark}
\newtheorem{solution}[theorem]{Solution}
\newtheorem{summary}[theorem]{Summary}

\begin{frontmatter}

\title{Stochastic comparisons of stratified sampling techniques for some Monte Carlo estimators}
\runtitle{Stratified sampling for some Monte Carlo estimators}

\begin{aug}
\author[1]{\fnms{Larry} \snm{Goldstein}\thanksref{1}\ead[label=e1]{larry@math.usc.edu}},
\author[2]{\fnms{Yosef} \snm{Rinott}\thanksref{2}\ead[label=e2]{rinott@mscc.huji.ac.il}}
\and
\author[3]{\fnms{Marco} \snm{Scarsini}\corref{}\thanksref{3}\ead[label=e3]{marco.scarsini@luiss.it}}

\runauthor{L. Goldstein, Y. Rinott and M. Scarsini}

\address[1]{Department of Mathematics,
University of Southern California,
Kaprielian Hall, Room 108,
3620 Vermont Avenue, Los Angeles,
CA 90089-2532,
USA.
\printead{e1}}

\address[2]{Department of Statistics and Center for the Study of Rationality,
Hebrew University of Jerusalem,
Mount Scopus,
Jerusalem 91905,
Israel and LUISS,
Roma,
Italy.\\
\printead{e2}}

\address[3]{Dipartimento di Scienze Economiche e Aziendali LUISS,
Viale Romania 12,
I--00197 Roma,
Italy and HEC,
Paris,
France.
\printead{e3}}
\end{aug}

\received{\smonth{4} \syear{2009}}
\revised{\smonth{5} \syear{2010}}

\begin{abstract}
We compare estimators of the (essential) supremum and the integral of a
function $f$ defined on a measurable space when $f$ may be observed at
a sample of points in its domain, possibly with error. The estimators
compared vary in their levels of stratification of the domain, with the
result that more refined stratification is better with respect to
different criteria. The emphasis is on criteria related to stochastic
orders. For example, rather than compare estimators of the integral of
$f$ by their variances (for unbiased estimators), or mean square error,
we attempt the stronger comparison of convex order when possible. For
the supremum, the criterion is based on the stochastic order of
estimators.
\end{abstract}

\begin{keyword}
\kwd{convex loss}
\kwd{convex order}
\kwd{majorization}
\kwd{stochastic order}
\kwd{stratified sampling}
\end{keyword}

\end{frontmatter}

\section{Introduction}\label{Se:Int}

In many situations, the cost of computing the value of a function $f$
is very high, because either the analytic expression of the function is
extremely complex or the value is the result of a costly experiment.
For example, $f$ could be the level of toxicity as a reaction to
different doses of certain drugs, the output of a chemical experiment,
or the  survival time of a patient undergoing a certain treatment.
Therefore the function can be computed only at a limited number of
points. One standard way to choose these points is via some Monte Carlo
randomization. Different possibilities arise: points could be sampled
totally at random or some stratification could be used. When properly
carried out, stratification is known to improve the performance of
estimators. The purpose of this paper is to qualify the above statement
in some relevant cases and compare different sampling stratifications
according to some suitable criteria.

Often the object of interest is some functional of $f,$ such as its
supremum or integral. Monte Carlo estimation of such functionals  is
the subject of a very large number of papers. In most cases some
regularity of the function $f$ is assumed; see, for example,
\cite{NovdsdLNMSpringer1988,ZhiChedsMCM1996}. Under some regularity
conditions it is often reasonable to estimate the entire function and
then use a plug-in method to estimate the functional. When no
regularity is assumed for $f$, then it may be more reasonable to
estimate the functional directly.

Given a measurable space $(\mathfrak{U}, \mathcal{U})$, let
$f\dvtx\mathfrak{U} \to \mathbb{R}$ be a measurable function $f$. In
order to estimate $\theta :=\sup_{x\in\mathfrak{U}}f(x),$ we can draw a
sample $X_1, \dots, X_n$ of $n$ points in $\mathfrak{U}$ and use the
estimator  $T :=\max(f(X_1), \ldots, f(X_n))$.  Alternatively we can
sample the $X$'s by resorting to some stratification. Ermakov,
Zhiglyavski{\u\i} and Kondratovich \cite{ErmZhiKondsdDANSSSR1988},
Kondratovich and Zhigljavsky \cite{KonZhidsdLNS1998} and Zhigljavsky and
{\v{Z}}ilinskas \cite{ZhiZilfdSPRINGER2008} prove that, if we consider
two partitions of $\mathfrak{U}$, one of which is a refinement of the
other, and we sample in proportion to the measure of each element of
the partition, then the more refined partition produces a
stochastically larger estimator of the supremum. Since these estimators
are  almost surely smaller than $\theta$ (hence biased) and consistent,
the stochastically larger one performs better. Thus, the more we
stratify, the better the estimator we obtain.

In our paper we extend this result and show that the stochastic
comparison for estimators of the supremum holds also when observations
are censored, that is, when for a sample of pairs of random variables
$(U_i, Z_i)$ we only know whether $Z_i \le f(U_i)$ or not. In
applications, there may be situations where exact evaluation of $f(u)$
at a given point is difficult or expensive, whereas a comparison of
$f(u)$ to a given constant $t$ is (at least for most values of $t$)
much easier.  For example, if $f(u)$ represents a lifetime, it may be
easier to see if it has exceeded a certain value, rather than wait  to
obtain the exact value $f(u)$ itself. This amounts to censoring.

When we want to estimate the integral  $I(f)$ of the function $f$, then
it is easy to construct an unbiased estimator of $I(f)$ by using
different stratified samples. Unbiasedness of these estimators implies
that the comparison criterion cannot be the stochastic order, as used
for the maximum.

In much of the literature estimators are compared in terms of a given
loss function, which may be arbitrary. Typically the loss function is
quadratic, so the criterion is the mean square error, that is, the
variance, when the estimator is unbiased. More generally, it may be
possible to find comparison criteria that are valid for large classes
of loss functions; for instance, all losses of the type $|W-I(f)|^{p}$,
where  $W$ is an estimator of $I(f)$ and $p \ge 1$, or even the class
of all convex loss functions. The use of the entire class of convex
loss functions in inference goes back at least to \cite{LaydsdJRSSB1972}
and \cite{LaySildsdB1968}. Similar ideas were later used  by Berger
\cite{BerdsdAS1976}, Kozek \cite{KozdsdJMVA1977}, Lin and Mousa
\cite{LinMoudsdAISM1982}, Eberl \cite{EbedsdSD1984}, Bai and Durairajan
\cite{BaiDurdsdJSPI1997}, and Petropoulos and
Kourouklis \cite{PetKoureAISM2001}. A comparison of the performance of
different estimators, with respect to all convex loss functions, can be
achieved by considering the convex order. Comparison of experiments in
terms of the convex order traces back to
\cite{BladsdP2BSMSP1951,BladsdAMS1953}.

It is well known that stratification reduces the variance of estimators
of $I(f)$, but, as will be shown below,  stratification does not
necessarily reduce $\mathbb{E}[|W-I(f)|^{p}]$, for $p \neq 2$, which
implies that, even if stratification is useful in $L_{2}$, it may be
counterproductive in $L_{1}$. We will show that in some circumstances
stratified sampling is better not just in $L_{2}$, but in terms of the
convex order, which in turn implies that it is better in $L_{p}$ for
every $p \ge 1$. This is the case when observations are censored, the
function $f$ is univariate and monotone, or the function is
multivariate and monotone and the sampling is independent across
coordinates. Papageorgiou \cite{PapdsdJC1993} shows the computational
advantage of using randomized methods to compute the integral of
monotone $d$-variate functions, and shows how this depends on $d$.

Our results also hold when the function $f$ can only be observed with
noise; for instance, when $f$ is observed as the outcome of some
experiment. Moreover, our regularity assumptions on the function $f$
are rather non-restrictive: measurability when estimating the maximum,
boundedness when observations are censored, and sometimes monotonicity
when estimating the integral.

We emphasize that, in our framework, evaluation of $f$ by experiment is
the costly part and any precalculations, such as those required for
computing strata and sampling from the conditional distributions in
strata, even if computer-time consuming, are considered to have a
relatively negligible cost.

The paper is organized as follows. Section~\ref{Se:notation} fixes
notation and reviews various properties of stochastic orders and
certain dependence structures. Section~\ref{se:supcensor} compares
estimators of the supremum of a function, considering also the case of
censored observations. Section~\ref{se:integralcensor} compares
estimators of integrals: First a variance comparison is shown to hold
in general, even when observations are affected by errors. Then a
counterexample is provided for a non-quadratic loss function.  Then
censored observations are considered and a comparison in terms of the
convex order is proved in this case. Finally, monotone functions are
examined. In the univariate case, a convex order comparison holds. In
the multivariate case, this is true under some additional conditions on
the stratification and on the dependence of the underlying random
vector.

Numerical examples can be found in \cite{GolRinScadsdmimeo2010}.

\section{Notation and preliminaries}\label{Se:notation}

In this paper a probability space $(\Omega, \mathcal{F}, \mathbb{P})$
is assumed in the background. The \textit{stochastic order}
$\le_{\st}$, the \textit{convex order} $\le_{\cx}$, the
\textit{increasing convex order} $\le_{\icx}$, and the
\textit{majorization order} $\prec$ are defined as follows (see, e.g.,
\cite{MarOlkdsdACADEMIC1979,MueStodsdWiley2002,ShaShadsfSPRI2007}). Given two random vectors $\mathbf{X},
\mathbf{Y}$, we say that $\mathbf{Y} \le_{\st} \mathbf{X}$ if
\begin{equation}\label{eq:EphiXEphiY}
\mathbb{E}[\phi(\mathbf{Y})] \le \mathbb{E}[\phi(\mathbf{X})]
\end{equation}
for all non-decreasing functions $\phi$.  We say that $\mathbf{Y}
\le_{\cx} \mathbf{X}$ if \eqref{eq:EphiXEphiY} holds for all convex
functions $\phi$ and $\mathbf{Y} \le_{\icx} \mathbf{X}$ if
\eqref{eq:EphiXEphiY} holds for all non-decreasing convex functions
$\phi$. It is well known that  $\mathbf{Y} \le_{\st} \mathbf{X}$ iff
$\mathbb{P}(\mathbf{Y} \in A) \le \mathbb{P}(\mathbf{X} \in A)$ for all
increasing sets $A$, where we call a set \textit{increasing} if its
indicator function is non-decreasing. In the case of univariate random
variables $X, Y$, the above inequality becomes $\mathbb{P}(Y \le t) \ge
\mathbb{P}(X \le t)$ for all $t \in \mathbb{R}$. It is well known that
$X \le_{\cx} Y$ implies $\mathbb{E}[X] = \mathbb{E}[Y]$ and $\Var[X]
\le \Var[Y]$.

The statement $\mathbf{Y} \le_{\st} \mathbf{X}$ depends only on the
marginal laws $\mathcal{L}(\mathbf{Y})$ and $\mathcal{L}(\mathbf{X})$,
so sometimes we write $\mathcal{L}(\mathbf{Y}) \le_{\st}
\mathcal{L}(\mathbf{X})$, and analogously for $\le_{\cx}$ and
$\le_{\icx}$.

Given two vectors $\mathbf{x} = (x_{1}, \ldots, x_{n})$, $\mathbf{y} =
(y_{1}, \ldots, y_{n})$, we write $\mathbf{y} \prec \mathbf{x}$ if
\[
\sum_{i=1}^{k} y_{i}^\downarrow \le \sum_{i=1}^{k}
x_{i}^\downarrow\qquad \mbox{for } k=1,\ldots, n-1,
\qquad\sum_{i=1}^{n} y_{i} =\sum_{i=1}^{n} x_{i},
\]
where $y_{1}^\downarrow \ge \cdots \ge y_{n}^\downarrow$ is the
decreasing rearrangement of $\mathbf{y}$, and analogously for
$\mathbf{x}$. The relation $\mathbf{y} \prec \mathbf{x}$ holds if and
only if there exists an $n \times n$ doubly stochastic matrix
$\mathbf{D}$ such that $\mathbf{y} = \mathbf{D} \mathbf{x}$.

A function $\psi\dvtx \mathbb{R}^n \rightarrow \mathbb{R}$ is called
Schur convex or Schur concave if $\mathbf{y} \prec \mathbf{x}$ implies
$\psi(\mathbf{y}) \le \psi(\mathbf{x})$ or $\psi(\mathbf{y}) \ge
\psi(\mathbf{x})$, respectively. If $\varphi\dvtx \mathbb{R}
\rightarrow \mathbb{R}$ is convex then $\psi(\mathbf{x})=\sum_{i=1}^n
\varphi(x_i)$ is Schur convex.

A random vector $\mathbf{X}$ is \textit{associated} if for all
non-decreasing functions $\phi, \psi$ we have
$\Cov[\phi(\mathbf{X}),\psi(\mathbf{X})] \ge 0$.

Recall that a subset $A \subset \mathbb{R}^{d}$ is a \textit{lattice}
if it is closed under componentwise maximum $\vee$ and minimum
$\wedge$. A random vector $\mathbf{X}$ is \textit{multivariate totally
positive of order $2$} (MTP$_2$) if its support is a lattice and its
density $f_{\mathbf{X}}$ with respect to some product measure on
$\mathbb{R}^{d}$ satisfies $f_{\mathbf{X}}(\mathbf{s})
f_{\mathbf{X}}(\mathbf{t}) \le f_{\mathbf{X}}(\mathbf{s} \vee
\mathbf{t}) f_{\mathbf{X}}(\mathbf{s} \wedge \mathbf{t})$ for all
$\mathbf{s}, \mathbf{t} \in \mathbb{R}^d$. MTP$_{2}$ implies
association. Also, any vector having independent components is
MTP$_{2}$.

Let $U$ be a random variable with values in some measurable space
$(\mathfrak{U}, \mathcal{U})$ with non-atomic law $P_{U}$. A finite
sequence $\mathcal{B} = (B_{1}, \dots, B_{b})$ of subsets of
$\mathfrak{U}$ is called an \textit{ordered partition} of
$\mathfrak{U}$ if $B_{i} \cap B_{j} = \varnothing$ for $i, j \in \{1,
\dots, b\}$, $i \neq j$, and $\bigcup_{i=1}^{b} B_{i} = \mathfrak{U}$.
For the sake of brevity in the sequel, whenever we say  ``partition''
we mean ``ordered partition.''

Here we consider partitions $\mathcal{B} = (B_{1}, \dots, B_{b})$ of
$\mathfrak{U}$, where the sets  $B_{i}$ are measurable and such that
for $i=1, \dots, b$ we have $\mathbb{P}(U \in B_{i})= k_{i}/n$ for some
$k_{i} \in \{1,\ldots,n\}$ satisfying $\sum_i k_i = n$. We say that
such a partition $\mathcal{B}$ of $\mathfrak{U}$ and a partition
$\mathcal{B}^{*}= (B_{1}^{*}, \dots, B_{b}^{*})$ of $N :=
\{1,\ldots,n\}$ are associated if the cardinalities $|B_i^{*}|$ of the
sets $B_i^{*}$ satisfy $|B_i^{*}|=k_i$ for $i=1, \dots, b$. We then
have
\begin{equation}\label{eq:PBB*}
\mathbb{P}(U \in B_{i}) = \frac{|B_{i}^{*}|}{n}.
\end{equation}
The notation $B \in \mathcal{B}$ means that $B$ is one of the sets
$B_{i}$ that comprise $\mathcal{B}$ and, given $B \in \mathcal{B}$, we
let $B^*$ denote the corresponding set $B_i^{*}$ in $\mathcal{B}^{*}$
such that \eqref{eq:PBB*} holds.

Given two partitions $\mathcal{B}^{*}= (B_{1}^{*}, \dots, B_{b}^{*})$
and $\mathcal{C}^{*} = (C_{1}^{*}, \ldots, C_{c}^{*})$ of $N$, we write
$\mathcal{C}^{*} \le_{\refin} \mathcal{B}^{*}$; that is, that
$\mathcal{B}^{*}$ is a refinement of $\mathcal{C}^{*}$ when every set
in $\mathcal{C}^{*}$ is the union of sets in $\mathcal{B}^{*}$. We will
use the same order $\le_{\refin}$  for partitions of $\mathfrak{U}$.
Clearly, if $\mathcal{C}$ and $\mathcal{B}$ are partitions of
$\mathfrak{U}$, each of which can be associated to some partition of
$N$, then $\mathcal{C} \le_{\refin} \mathcal{B}$ implies that there
exist partitions $\mathcal{C}^{*}$ and $\mathcal{B}^{*}$ associated to
$\mathcal{C}$ and $\mathcal{B}$, respectively, satisfying
$\mathcal{C}^{*} \le_{\refin} \mathcal{B}^{*}$.

Call $\mathcal{A}^{*} = (\{1\}, \dots, \{n\})$ the finest partition of
$N$ and $\mathcal{D}^{*}=(N)$ the coarsest partition  of $N$. Then
$\mathcal{D}^{*}\le_{\refin} \mathcal{B}^{*}\le_{\refin}
\mathcal{A}^{*}$ for all $\mathcal{B}^{*}$, and for any partition
$\mathcal{A}$ of $\mathfrak{U}$ associated to $\mathcal{A}^{*}$ we have
$\mathbb{P}(U \in A_{i}) = 1/n$.

For a partition $\mathcal{B}$ and $B \in \mathcal{B}$, let $P_{U|B}$
denote the conditional law of $U$ given $U\in B$. Let $\{V_{j}^{B}, j
\in B^* \}$ be random variables with law $P_{U|B}$ with $\{V_{j}^{B}, j
\in B^*, B \in \mathcal{B}\}$ independent.

\section{The supremum}\label{se:supcensor}

Let $f\dvtx \mathfrak{U} \to \mathbb{R}$ be measurable, and define
\begin{equation}
\label{eq:WS} W_{\mathrm{S}}^{\mathcal{B}} = \max\limits_{B \in
\mathcal{B}}\max\limits_{j \in B^{*}} f(V_{j}^{B}),
\end{equation}
where the subscript $\mathrm{S}$ indicates that
$W_{\mathrm{S}}^{\mathcal{B}}$ will be used to estimate the (essential)
supremum of the function $f$.

Given a random variable $U$  with values in  $(\mathfrak{U},
\mathcal{U})$, let $f^{*} := \operatorname{ess\,sup}f(U)$. It is clear
that for any choice of partition $\mathcal{B}$,
$\mathbb{P}(W_{\mathrm{S}}^{\mathcal{B}} \le f^{*})=1$. The following
result compares two estimators of type $W_{\mathrm{S}}^{\mathcal{B}}$.
Since both estimators underestimate $f^{*}$, the stochastically larger
one is preferable. This theorem, which goes back to
\cite{ErmZhiKondsdDANSSSR1988} and \cite{KonZhidsdLNS1998}, can also be
found in \cite{ZhiZilfdSPRINGER2008}, Theorem~3.4.

\begin{theorem}\label{th:maxgeneral}
If $\mathcal{C} \le_{\refin}\mathcal{B}$, then
$W_{\mathrm{S}}^{\mathcal{C}} \le_{\st} W_{\mathrm{S}}^{\mathcal{B}}$.
\end{theorem}

A short proof of Theorem~\ref{th:maxgeneral}, different from the one in
the
\cite{ZhiZilfdSPRINGER2008}, can be found in the
\hyperref[se:appendix]{Appendix}.

As mentioned in the Section \ref{Se:Int}, data are not always observed exactly in
many practical situations, but may be censored for various reasons,
including budget constraints. We extend now the comparison result of
Theorem~\ref{th:maxgeneral} to the case of censored observations. Let
$f\dvtx \mathfrak{U} \to \mathbb{R}$ be  bounded; without loss of
generality, we take $0 \le f(u) \le 1$ for all $u \in \mathfrak{U}$. In
this section we assume that, for a sample of points of the type $(u,t)
\in \mathfrak{U} \times [0,1]$, we are allowed to observe only the
value of $t$ and whether $t > f(u)$.

For any partition  $\mathcal{B}$ with associated partition
$\mathcal{B}^*$,  let $\{V_{j}^{B}, j \in B^{*}\}$, $B \in \mathcal{B}$
and $\{T_{j}, j \in N\}$ be independent random variables with law
$P_{U|B}$ and the uniform distribution on $[0,1]$, respectively, and
let
\[
S^{\mathcal{B}} = \bigcup_{B\in \mathcal{B}}   \{j \in B^{*} \dvtx
T_{j} \le f(V_{j}^{B})\}\quad  \mbox{and}\quad  W_{\CS}^{\mathcal{B}} =
\max\limits_{j \in S^{\mathcal{B}}} T_{j}.
\]
When $S^{\mathcal{B}} = \varnothing$ we set $W_{\CS}^{\mathcal{B}} =
0$. The letter C in the subscript $\CS$ indicates censored data. It is
clear that $\mathbb{P}(W_{\CS}^{\mathcal{B}}  \le f^{*}) = 1$, so the
estimator $W_{\CS}^{\mathcal{B}}$  underestimates $f^{*}$.

\begin{theorem}\label{th:imax}
If $\mathcal{C} \le_{\refin} \mathcal{B}$, then $W_{\CS}^{\mathcal{C}}
\le_{\st} W_{\CS}^{\mathcal{B}}$.
\end{theorem}

\begin{pf}
Below, when we write $V_{j}^{B}$ without specifying $B$, we
mean that $B \in \mathcal{B}$ corresponds in the sense of
\eqref{eq:PBB*} to  the set $B^{*} \in \mathcal{B}^{*}$,  which
contains the index $j$. For any $t \in [0,1],$ we may calculate the
distribution function of $ W_{\CS}^{\mathcal{B}}$ at $t$ by writing
\begin{eqnarray*}
\{W_{\CS}^{\mathcal{B}} \le t\}
&=&
\bigcup_{R \subset N}\Bigl\{\max_{j \in S^{\mathcal{B}}} T_{j} \le t, S^{\mathcal{B}}=R\Bigr\}\\
&=&
\bigcup_{R \subset N} \{T_{j} \le t, T_{j} \le f(V_{j}^{B})\mbox{ for all $j \in R$, and } T_{j} > f(V_{j}^{B}) \mbox{ for all $j \notin R$} \} \\
&=&
\bigcup_{R \subset N} \{T_{j} \le t \wedge f(V_{j}^{B}) \mbox{ for all $j \in R$, and } T_{j} > f(V_{j}^{B}) \mbox{ for all $j \notin R$}\}.
\end{eqnarray*}
Hence, conditionally on $\{V_{j}^{B}$, $j \in B^{*}$, $B \in
\mathcal{B}\}$, using the fact that the $T_{j}$'s are uniform, we
obtain:
\begin{eqnarray}\label{eq:Z-cond-dist}
&&\mathbb{P}(W_{\CS}^{\mathcal{B}} \le t \vert V_{j}^{B}, j \in B^{*}, B \in\mathcal{B})\nonumber
\\
&&\quad=
\sum_{R \subset N}\prod_{j\in R}\mathbb{P}\bigl(T_{j}\le t\wedge f(V_{j}^{B})\bigr)\prod_{j\notin R}\mathbb{P}\bigl(T_{j}>f(V_{j}^{B})\bigr)\nonumber
\\[-8pt]\\[-8pt]
&&\quad=
\sum_{R\subset N}\prod_{j\in R}\bigl(t\wedge f(V_{j}^{B})\bigr)\prod_{j\notin R}\bigl(1- f(V_{j}^{B})\bigr)\nonumber
\\
&&\quad=
\sum_{h_{1}=1}^{|B_{1}^{*}|}\dots\sum_{h_{b}=1}^{|B_{b}^{*}|}\mathop{\sum_{R \subset N}}_{\forall i,|R\cap B_i^{*}|=h_i}
\prod_{j \in R}\bigl(t\wedge f(V_{j}^{B})\bigr)\prod_{j\notin R}\bigl(1- f(V_{j}^{B})\bigr).\nonumber
\end{eqnarray}
Taking expectation we obtain the unconditional distribution,
\begin{eqnarray*}
\mathbb{P}(W_{\CS}^{\mathcal{B}}\le t)
&=&
\sum_{h_{1}=1}^{|B_{1}^{*}|}\cdots\sum_{h_{b}=1}^{|B_{b}^{*}|}\prod_{i=1}^{b}\pmatrix{|B_{i}^{*}|\cr h_{i}}\biggl(\int_{B_{i}}\bigl(t \wedge f(u)\bigr)\,\diff P_{U|{B_{i}}}(u)\biggr)^{h_{i}}
\\
&&\hphantom{\sum_{h_{1}=1}^{|B_{1}^{*}|}\cdots\sum_{h_{b}=1}^{|B_{b}^{*}|}\prod_{i=1}^{b}}{}\times\biggl(\int_{B_{i}}\bigl(1-f(u)\bigr)\,\diff P_{U|{B_{i}}}(u)\biggr)^{|B_{i}^{*}|-h_{i}}
\\
&=& \prod_{B \in \mathcal{B}}\biggl(\int_{B}\bigl(t \wedge f(u)\bigr)\,\diff P_{U|{B}}(u)+\int_{B}\bigl(1-f(u)\bigr)\,\diff P_{U|{B}}(u)\biggr)^{|B^{*}|}.
\end{eqnarray*}
Let
\begin{eqnarray*}
q^{B}
&=&
\int_{B}\bigl(t \wedge f(v)\bigr)\,\diff P_{U|{B}}(v)+\int_{B}\bigl(1-f(v)\bigr)\,\diff P_{U|{B}}(v)
\\
&=&
\int_{B_{}} \bigl[\bigl(t \wedge f(v)\bigr)+\bigl(1-f(v)\bigr)\bigr]\,\diff P_{U|{B}}(v).
\end{eqnarray*}
If $C$ is a union of disjoint sets $B_i$, then
\begin{equation}\label{eq:qCqB}
q^{C}=\sum_i q^{B_i}\frac{\mathbb{P}(U \in B_i)}{\mathbb{P}(U \in C)} =
\sum_i q^{B_i}\frac{|B_i^*|}{|C^*|}.
\end{equation}

If $\mathcal{C} \le_{\refin} \mathcal{B}$, then
\[
(\underbrace{q^{C_{1}},\ldots,q^{C_{1}}}_{|C_{1}^{*}|}, \ldots,
\underbrace{q^{C_{c}},\ldots,q^{C_{c}}}_{|C_{c}^{*}|}) \prec
(\underbrace{q^{B_{1}},\ldots,q^{B_{1}}}_{|B_{1}^{*}|}, \ldots,
\underbrace{q^{B_{b}},\ldots,q^{B_{b}}}_{|B_{b}^{*}|}).
\]
To see this, observe that \eqref{eq:qCqB} implies that the vector on
the left-hand side above is obtained from the one on the right by
multiplying it by the $n \times n$ doubly stochastic matrix
$\mathbf{D}$, which is block diagonal where the $i$th block is the
$|C_{i}^{*}| \times |C_{i}^{*}|$ matrix with all entries equal to
$1/|C_{i}^{*}|$. Therefore, by the Schur concavity of the function
$(\theta_{1}, \dots, \theta_{n}) \mapsto \prod_{i=1}^n \theta_i$, we
have
\[
\mathbb{P}(W_{\CS}^{\mathcal{C}} \le t) = \prod_{C \in \mathcal{C}}
(q^{C}  )^{|C^{*}|} \ge \prod_{B \in \mathcal{B}} (q^{B} )^{|B^{*}|} =
\mathbb{P}(W_{\CS}^{\mathcal{B}} \le t) .
\]
\upqed
\end{pf}

For every $n \in \mathbb{N}$ and for every partition $\mathcal{B}_n$
associated to a partition $\mathcal{B}_n^{*}$ of $\{1, \ldots, n\}$, we
have $W_{\CS}^{\mathcal{B}_n} \le_{\st}
W_{\mathrm{S}}^{\mathcal{B}_n}$. Therefore,
\[
W_{\CS}^{\mathcal{D}_n} \le_{\st} W_{\CS}^{\mathcal{B}_n}  \le_{\st}
W_{\mathrm{S}}^{\mathcal{B}_n}  \le_{\st} f^{*}.
\]
Since $W_{\CS}^{\mathcal{D}_n}$ is consistent for $f^{*}$ as $n
\rightarrow \infty$, we have that $W_{\CS}^{\mathcal{B}_n}$ and
$W_{\mathrm{S}}^{\mathcal{B}_n}$ are consistent, too.

\section{The integral}\label{se:integralcensor}

With the subscript $\integr$ standing for integral, let
\begin{eqnarray}
\label{eq:WI}W_{\integr}^{\mathcal{B}}
&=&
\frac{1}{n} \sum_{B \in \mathcal{B}}\sum_{j \in B^{*}} f(V_{j}^{B}),
\\
\label{eq:WIE}W_{\IE}^{\mathcal{B}}
&=&
\frac{1}{n} \sum_{B \in\mathcal{B}}\sum_{j \in B^{*}}\bigl(f(V_{j}^{B})+\varepsilon_{j}\bigr),
\end{eqnarray}
where  the variables $\varepsilon_{j}$ are independent copies of a
random variable $\varepsilon$ having mean $0$ and finite variance,
independent of the variables $V_{j}^{B}$. Clearly
$W_{\integr}^{\mathcal{B}}$ and $W_{\IE}^{\mathcal{B}}$ are both
unbiased estimators of $\overline{f} :=\mathbb{E}[f(U)]=\int f(U)\,\diff
\mathbb{P}$ when $\int |f(U)|\,\diff \mathbb{P}$  is finite, and
$W_{\integr}^{\mathcal{B}}$ is the special case of
$W_{\IE}^{\mathcal{B}}$ when the error has zero variance; that is,
there  is no measurement error.

The following result is well known when the error has zero variance
(see, e.g., \cite{GladsdSpringer2004}, Section~4.3). We extend it to a
more general case, relevant when the evaluation of $f$ is the result of
an experiment.

\begin{theorem}\label{th:interror}
If $\mathcal{C} \le_{\refin} \mathcal{B}$, then
$\Var[W_{\IE}^{\mathcal{B}}] \le \Var[W_{\IE}^{\mathcal{C}}]$.
\end{theorem}

The proof of Theorem~\ref{th:interror} can be found in the
\hyperref[se:appendix]{Appendix}.

It follows immediately from Theorem~\ref{th:interror} that
$\Var[W_{\IE}^{\mathcal{A}}] \le \Var[W_{\IE}^{\mathcal{D}}]$, hence,
in particular, $\Var[W_{\integr}^{\mathcal{A}}] \le
\Var[W_{\integr}^{\mathcal{D}}]$. The following counterexample shows,
nevertheless, that, even when the function is observed without error,
$W_{\integr}^{\mathcal{A}} \not \le_{\cx} W_{\integr}^{\mathcal{D}}$;
that is,  domination in the convex order does not hold. In the
counterexample we consider the absolute error, that is, ($L_1$), rather
than mean square error, ($L_2$).

\begin{example}\label{counterexample}
Let $\mathfrak{U} = [0,1]$ and $U$ have a uniform distribution on
$[0,1]$. Furthermore, let   $n=2$,  $A_1=[0,1/2], A_2=(1/2,1]$. Define
\[
f(u) = 4I_{[0,1/2]}(u) + 2I_{(1/2,3/4]}(u) + 6I_{(3/4,1]}(u).
\]

Then $W_{\integr}^{\mathcal{D}}$  takes the values $2,3,4,5,6$ with
probabilities $(1,4,6,4,1)/16$, respectively. The variable
$W_{\integr}^{\mathcal{A}}$, based on one random observation from each
of the above intervals $A_i$, takes the values 3 and 5 each with
probability $1/2$. Therefore,
$\mathbb{E}[W_{\integr}^{\mathcal{A}}]=4=\mathbb{E}[W_{\integr}^{\mathcal{D}}]$.

We have $\Var[W_{\integr}^{\mathcal{D}}] =
\Var[W_{\integr}^{\mathcal{A}}] = 1$, but for the convex function
$\psi(u) = |u-4|$ we have
\[
\mathbb{E}[\psi(W_{\integr}^{\mathcal{D}})]
=\mathbb{E}|W_{\integr}^{\mathcal{D}}-4| = 2 \frac{2}{16}+2
\frac{4}{16}=\frac{12}{16} < 1 =
\mathbb{E}|W_{\integr}^{\mathcal{A}}-4| =
\mathbb{E}[\psi(W_{\integr}^{\mathcal{A}})].
\]

A more general example can be constructed as follows. Consider a
partition $\mathcal{A}$ associated to the finest partition
$\mathcal{A}^{*}$ of $N$. Split $A_{1}$ into two measurable subsets
$A_{1a}, A_{1b}$ such that $\mathbb{P}(U \in A_{1a}) = \mathbb{P}(U \in
A_{1b}) = 1/(2n)$. Consider now a function $f$ defined as follows:
\begin{equation}\label{eq:nonmonotonefunction}
f(u) =
\cases{
1, &\quad if $u \in A_{1a}$,\cr
-1, &\quad if $u \in A_{1b}$,\cr
0, &\quad elsewhere.
}
\end{equation}
For all $i \in N$ we have $\mathbb{E}[f(U)\mid U\in A_{i}] = 0$ and
\[
\Var[f(U)\mid U\in A_{i}] =
\cases{
1, &\quad for $i=1$,\cr
0, &\quad for $i \neq 1$.
}
\]
Hence
\[
\Var[W_{\integr}^{\mathcal{A}}]=\mathbb{E}[(W_{\integr}^{\mathcal{A}})^{2}] = \frac{1}{n^{2}}.
\]
Moreover, if $V_{1}, \dots, V_{n}$ are i.i.d. copies of $U$,
\[
\Var[W_{\integr}^{\mathcal{D}}]
=
\Var\Biggl[\frac{1}{n} \sum_{j=1}^{n} f(V_{j})\Biggr]
=
\frac{1}{n^{2}} \sum_{j=1}^{n} \Var[f(V_{j})]
=
\frac {1}{n^{2}}
=
\Var[W_{\integr}^{\mathcal{A}}].
\]

Analogously
\[
\mathbb{E}[|f(U)|\mid U\in A_{i}] =
\cases{
1, &\quad for $i=1$,\cr
0, &\quad for $i \neq 1$.
}
\]
Therefore
\[
\mathbb{E}|W_{\integr}^{\mathcal{A}}|
=
\sqrt{\mathbb{E}[(W_{\integr}^{\mathcal{A}})^{2}]}=\frac{1}{n}.
\]
For any square integrable random variable $Y$ we have $\mathbb{E}|Y|
\le \sqrt{\mathbb{E}[Y^{2}]}$ and the inequality is strict if $Y$ is
not almost surely constant. Hence
\[
\mathbb{E}|W_{\integr}^{\mathcal{D}}| <
\sqrt{\mathbb{E}[(W_{\integr}^{\mathcal{D}})^{2}]} =
\sqrt{\mathbb{E}[(W_{\integr}^{\mathcal{A}})^{2}]} =
\mathbb{E}|W_{\integr}^{\mathcal{A}}|  =\frac{1}{n}.
\]
\end{example}

Example~\ref{counterexample} proves that the convex order does not hold
in general between estimators $W_{\integr}^{\mathcal{B}}$ and
$W_{\integr}^{\mathcal{C}}$ when $\mathcal{C} \le_{\refin}
\mathcal{B}$. Nevertheless, in the following subsections we show that
under some natural conditions comparisons in the convex order are
possible.

\subsection{Censored observations}

Keeping the notation and spirit of  Section~\ref{se:supcensor},
consider a function $f$ such that $0 \le f(u) \le 1$ for all $u \in
\mathfrak{U}$. Assume  that for a sample of points of the type $(u,t)
\in \mathfrak{U} \times [0,1]$ we are allowed to observe only the value
of $t$ and whether $t \le f(u).$ Let
\[
W_{\CI}^{\mathcal{B}}=\frac{1}{n}\sum_{B\in\mathcal{B}}\sum_{j\in B^{*}}I_{\{T_{j}\le f(V_{j}^{B})\}}.
\]
Note that $W_{\CI}^{\mathcal{B}}$ is an unbiased estimator of
$\overline{f}=\mathbb{E}[f(U)]$, as
\begin{eqnarray*}
\mathbb{E}[W_{\CI}^{\mathcal{B}}]
&=&
\frac{1}{n} \sum_{B \in\mathcal{B}}\sum_{j \in B^{*}}\mathbb{P}\bigl(T_{j}\le f(V_{j}^{B})\bigr)
=
\frac{1}{n} \sum_{B \in \mathcal{B}}\sum_{j \in B^{*}}\int_{\mathfrak{U}} \int_0^1 I_{\{t \le f(u)\}}\,\diff t\,\diff P_{U|B}(u)
\\
&=&
\sum_{B \in \mathcal{B}}\frac{|B^{*}|}{n} \int_{\mathfrak{U}} f(u)\,\diff P_{U|B}(u)=\sum_{B \in \mathcal{B}}\mathbb{P}(B)\mathbb{E}[f(U)\mid U\in B]
\\
&=&
\mathbb{E}[f(U)].
\end{eqnarray*}

\begin{theorem}\label{th:censint}
If $\mathcal{C} \le_{\refin} \mathcal{B}$, then $W_{\CI}^{\mathcal{B}}
\le_{\cx} W_{\CI}^{\mathcal{C}}$.
\end{theorem}

\begin{pf}
By a result in \cite{KarNovdsdPJM1963} (see also
\cite{MarOlkdsdACADEMIC1979}, Sections 12.F and 15.E)  if
\[
X_{\mathbf{p}}= \frac{1}{n}\sum_{i=1}^n \xi_i,
\]
where $\xi_1,\ldots,\xi_n$ are  independent Bernoulli variables with
parameters  $p_1,\ldots,p_n$, and $\mathbf{p}=(p_1,\ldots,p_n)$, then
\begin{equation}\label{eq:KarNov}
\mathbf{p} \prec \mathbf{q}\quad\mbox{implies}\quad X_{\mathbf{q}} \le_{\cx} X_{\mathbf{p}}.
\end{equation}
Define
\[
p^{C}=\mathbb{P}\bigl(T_{j} \le f(V_{j}^{C})\bigr),\qquad p^{B}=\mathbb{P}\bigl(T_{j} \le f(V_{j}^{B})\bigr),
\]
and
\[
\mathbf{p}=(\underbrace{p^{C_{1}},\ldots, p^{C_{1}}}_{|C_{1}^{*}|},
\dots, \underbrace{p^{C_{c}},\ldots, p^{C_{c}}}_{|C_{c}^{*}|}),
\qquad
\mathbf{q}=(\underbrace{p^{B_{1}},\ldots, p^{B_{1}}}_{|B_{1}^{*}|},
\dots, \underbrace{p^{B_{b}},\ldots, p^{B_{b}}}_{|B_{b}^{*}|}).
\]
If $C = \bigcup_i B_i$, then
\[
p^{C}=\sum_i p^{B_i}\frac{|B_i|}{|C|},
\]
so $\mathbf{p} \prec \mathbf{q}$ and invoking \eqref{eq:KarNov}
completes the proof.
\end{pf}

Notice that in the case of censored observations, the comparison holds
in the convex order, whereas in the case of perfect observation, a
variance comparison holds, but Example~\ref{counterexample} shows that
comparisons in the convex order do not.

\subsection{Univariate monotone functions}\label{suse:uni}
In the rest of this subsection the space $\mathfrak{U}$ is totally
ordered and, without loss of generality, we choose $\mathfrak{U} =
[0,1]$. For subsets $G$ and $H$ of the real line, we write $G \le H$ if
$g \le h$ for every $g \in G$ and  $h \in H$. We call a partition
$\mathcal{B} = (B_{1}, \dots, B_{b})$ of $\mathfrak{U}$ monotone if
$B_{1}\le \cdots \le B_{b}$.

\begin{theorem}\label{th:increasinginterror}
Let $\mathcal{B}$ and $\mathcal{C}$ be   monotone partitions of
$\mathfrak{U}$ and let $\mathcal{C} \le_{\refin} \mathcal{B}$. If $f$
is non-decreasing, then
\begin{equation}\label{eq:increasinginterror}
W_{\IE}^{\mathcal{B}} \le_{\cx} W_{\IE}^{\mathcal{C}}.
\end{equation}
\end{theorem}

To prove Theorem~\ref{th:increasinginterror} we will apply the
following lemma.

\begin{lemma}\label{le:partmon1}
Let $\xi$ and $\eta$ be random variables such that $\xi \le_{\st}
\eta$, and  let $\xi_i$ and $\eta_{j}$ be independent copies of $\xi$
and $\eta$, respectively. Let $K$ be an integer-valued random variable,
independent of all $\xi_{j}$ and $\eta_{j}$, satisfying $K \le m$ for
some integer $m$ and having an integer-valued expectation,
$\mathbb{E}[K] = k$. Then
\begin{equation}\label{eq:stepart}
\sum_{j=1}^k \xi_{j} + \sum_{j=k+1}^{m} \eta_{j}  \le_{\cx}
\sum_{j=1}^K \xi_{j} + \sum_{j=K+1}^{m} \eta_{j}.
\end{equation}
\end{lemma}

\begin{pf} Since $\xi \le_{\st} \eta$ we may construct i.i.d. pairs $(\xi_{i},\eta_{i})$ with $\mathbb{P}(\xi_{i} \le \eta_{i})=1$ for all $i=1,\ldots,m$. We adopt the usual convention that if $k=0,$ then $\sum_{j=1}^k \xi_{j}=0$.
First note that, by Wald's lemma,
\[
\mathbb{E}\Biggl[\sum_{j=1}^k \xi_{j} + \sum_{j=k+1}^{m} \eta_{j}\Biggr]  =
\mathbb{E}\Biggl[ \sum_{j=1}^K \xi_{j} + \sum_{j=K+1}^{m} \eta_{j}\Biggr].
\]
Therefore (see, e.g., \cite{MueStodsdWiley2002}, Theorem~1.5.3) it
suffices to show that
\[
\sum_{j=1}^k \xi_{j} + \sum_{j=k+1}^{m} \eta_{j}  \le_{\icx}
\sum_{j=1}^K \xi_{j} + \sum_{j=K+1}^{m} \eta_{j}.
\]
Let $\phi$ be an increasing convex function and set
\[
g(k) := \mathbb{E}\Biggl[\phi\Biggl(\sum_{j=1}^k \xi_{j} +
\sum_{j=k+1}^{m}\eta_{j}\Biggr)\Biggr].
\]
Note that
\[
g(k)=\mathbb{E}\Biggl[\phi\Biggl(\sum_{j=1}^K \xi_{j}+\sum_{j=K+1}^{m}\eta_{j}\Biggr)\Bigr| K=k\Biggr]
\]
and
\[
\mathbb{E}[g(K)]=\mathbb{E}\Biggl[\phi\Biggl(\sum_{j=1}^K\xi_{j}+\sum_{j=K+1}^{m}\eta_{j}\Biggr)\Biggr].
\]
Thus we have to show that $g(k) \le\mathbb{E}[g(K)]$. Since
$\mathbb{E}[K]=k$, this follows readily by Jensen's inequality, once we
prove that $g(k)$ is a convex function.

The following part of the proof follows ideas of Ross and Schechner
\cite{RosSchfrdeOR1984}. Setting
\[
S_k=\sum_{j=1}^k \xi_{j} + \sum_{j=k+2}^{m} \eta_{j},
\]
we have
\[
g(k+1)-g(k) = \mathbb{E}[\phi(\xi_{k+1}+S_k)]-\mathbb{E}[\phi(\eta_{k+1}+S_k)].
\]
Since $\phi$ is convex, and $\xi_{k+1} \le \eta_{k+1}$,  the function
\[
h(s):=\mathbb{E}[\phi(\xi_{k+1}+S_k)\mid S_k=s]-\mathbb{E}[\phi(\eta_{k+1}+S_k)\mid S_k=s]
\]
is decreasing in $s$. Now note that
\[
S_{k+1}=\sum_{i=1}^{k+1} \xi_i + \sum_{i=k+3}^{m} \eta_i \le_{\st}S_k=\sum_{i=1}^k \xi_i + \sum_{i=k+2}^{m} \eta_i,
\]
because $\xi_{k+1} \le_{\st} \eta_{k+2}$. Hence $g(k+1)-g(k) =
\mathbb{E}[h(S_k)]$ is increasing in $k$, thus proving that $g$ is
convex, as required.
\end{pf}

\begin{pf*}{Proof of Theorem~\ref{th:increasinginterror}}
Since $\mathcal{B}=(B_1,\ldots,B_b)$ and $\mathcal{C}=(C_1,\ldots,C_c)$
are monotone partitions satisfying $\mathcal{C} \le_{\refin}
\mathcal{B}$, there exist $1=i_1<i_2<\cdots<i_c < i_{c+1}=b+1$ such
that
\[
C_q=\bigcup_{j=i_q}^{i_{q+1}-1}B_j\qquad\mbox{for }q=1,\ldots,c.
\]
As the union above may be formed by taking the union of two consecutive
sets at a time, it  suffices to prove \eqref{eq:increasinginterror} for
the case where $c=b-1$, $C_m=B_m \cup B_{m+1}$, $C_k=B_k$ for $k \in
\{1, \dots, m-1\}$, and $C_{k} = B_{k+1}$ for $k \in \{m+1, \dots,
c\}$.

In this case we have
\begin{eqnarray*}
W_{\IE}^{\mathcal{B}}
&=&
\frac{1}{n}\biggl[\sum_{C\neq C_{m}}\sum_{j\in C^{*}}f(V_{j}^{C})+\sum_{j\in B_{m}^{*}} f(V_{j}^{B_{m}})+\sum_{j\in B_{m+1}^{*}}f(V_{j}^{B_{m+1}})+\sum_{j\in N}\varepsilon_{j}\biggr],
\\
W_{\IE}^{\mathcal{C}}
&=&
\frac{1}{n}\biggl[\sum_{C \neq C_{m}} \sum_{j\in C^{*}} f(V_{j}^{C}) +  \sum_{j\in C_{m}^{*}} f(V_{j}^{C_{m}}) +\sum_{j\in N} \varepsilon_{j}\biggr].
\end{eqnarray*}
Note that
\[
\mathcal{L}\biggl(\sum_{j\in C_{m}^{*}} f(V_{j}^{C_{m}})\biggr) =
\mathcal{L}\Biggl(\sum_{j=1}^{K} f(V_{j}^{B_{m}}) +\sum_{j=K+1}^{|C_{m}^{*}|}
f(V_{j}^{B_{m+1}})\Biggr),
\]
where $K$ is binomially distributed with parameters
\[
\biggl(|C_{m}^{*}|,\frac{|B_{m}^{*}|}{|C_{m}^{*}|}\biggr).
\]
It is easy to see that if two variables are ordered by the convex order
(see \eqref{eq:EphiXEphiY}) and we add the same independent variable to
each one, to wit, $\sum_{j\in N} \varepsilon_{j}$, then the convex
order is preserved. This fact and Lemma~\ref{le:partmon1} now yield
\eqref{eq:increasinginterror}.
\end{pf*}

\subsection{Multivariate monotone functions}\label{suse:multi}
In this section we extend the results in Section \ref{suse:uni} to the
multivariate case. When we consider multivariate monotone functions,
stratifying can still yield  improvement in the convex order, but some
restrictions are needed, both on the distribution of the random vector
used for sampling and on the stratifying partitions. More specifically,
we consider estimation of an integral with respect to a random vector
whose components are independent and under a stratification that
preserves independence on each set of the partition. The result we
prove below actually only requires that the random vector have an
MTP$_{2}$ distribution (independence being a particular case of it) and
that the stratification  preserves MTP$_{2}$.

Let $f\dvtx[0,1]^d \rightarrow [0,1]$ be non-decreasing in each variable
and let $\mathbf{U}$ be a random vector taking values in $[0,1]^d$ with
a non-atomic distribution. Our goal is to show that the estimate of
$\mathbb{E}[f(\mathbf{U})]$ improves by refining stratifications as
follows. Recalling the definitions in Section \ref{Se:notation}, start
with a partition $\mathcal{C} = (C_{1}, \dots, C_{b})$ of $[0,1]^d$
such that for some $i$  the distribution $\mathcal{L}(\mathbf{U} \mid
\mathbf{U} \in C_{i})$ is associated. Then split $C_{i}$ into $C_{i}
\cap G$ and $C_{i} \cap G^{c}$, where $G$ is an increasing set.
Lemma~\ref{le:multconvexint} below shows that the new partition
obtained by this splitting achieves a better estimator of the integral
in terms of the convex order and Theorem~\ref{th:multconvexintgen}
provides some conditions for its application.

\begin{theorem}\label{th:multconvexintgen}
Consider a partition $\mathcal{C} = (C_{1}, \dots, C_{c})$ of
$[0,1]^{d}$ where each $C_i$ is a lattice. Let $\mathcal{B}$ be a
partition obtained by a sequence of refinements $\mathcal{C} =
\mathcal{C}_1 \le_{\refin} \cdots \le_{\refin} \mathcal{C}_m=
\mathcal{B}$, such that for $k=1,\ldots,m-1$ the partition
$\mathcal{C}_{k+1}$ is obtained from $\mathcal{C}_k$ by splitting one
set of $\mathcal{C}_k$, say $C_{i_k,k}$, into $C_{i_k,k} \cap G_k$ and
$C_{i_k,k}\cap G_k^{c}$, where $G_{k}= \{\mathbf{x}=(x_1,\ldots,x_d)
\in [0,1]^d\dvtx a_{k} \le x_j\} $ for some $a_{k} \in [0,1]$ and some $j
\in \{1,\ldots,d\}$.

If $\mathbf{U}$ is MTP$_{2}$ on $[0,1]^{d}$ and $f\dvtx[0,1]^d \rightarrow
[0,1]$ is non-decreasing, then $W_{\IE}^{\mathcal{B}}
\le_{\cx} W_{\IE}^{\mathcal{C}}$.
\end{theorem}

As mentioned earlier, independence is a particular (and in our
framework the most important) case of MTP$_2$. Independence makes
simulation of a multivariate random vector easy, even when conditioned
on an interval, since the strata can be constructed by knowing only the
quantiles of the marginal distributions. If the cost of simulation is
negligible relative to the cost of evaluating $f$, then even rejective
sampling can be used, once the strata are defined.

The proof of Theorem \ref{th:multconvexintgen} is preceded by the
following lemmas.

\begin{lemma}\label{le:assoc}
If $\mathbf{U}$ is an associated random vector, and $G$ is an
increasing set, then
\begin{equation} \label{eq:assocst}
\mathcal{L}(\mathbf{U} \mid\mathbf{U}\in G^c) \le_{\st}
\mathcal{L}(\mathbf{U} \mid \mathbf{U}\in G).
\end{equation}
Conversely, if \eqref{eq:assocst} holds for every increasing set $G$,
then $\mathbf{U}$ is associated.
\end{lemma}

\begin{pf}
First note that \eqref{eq:assocst} is equivalent to
\[
\mathbb{P}(\mathbf{U} \in A \vert \mathbf{U} \in G) \ge
\mathbb{P}(\mathbf{U}\in A \vert \mathbf{U} \in G^c)
\]
holding for all increasing sets $A$. The latter inequality is easily
seen to be equivalent to
\[
\mathbb{P}( \mathbf{U} \in A \cap G)[1-\mathbb{P}( \mathbf{U} \in G)]
\ge [\mathbb{P}(\mathbf{U} \in A )-\mathbb{P}( \mathbf{U} \in A \cap G)
]\mathbb{P}( \mathbf{U} \in  G).
\]
By simple cancelation this inequality is equivalent to
\[
\mathbb{P}(\mathbf{U} \in A \cap G) \ge \mathbb{P}(\mathbf{U} \in
A)\mathbb{P}(\mathbf{U} \in G),
\]
which is equivalent to association of the random vector $\mathbf{U}$
by, e.g., Shaked \cite{ShafdJMVA1982}.
\end{pf}

\begin{lemma}\label{le:multconvexint}
Consider a partition $\mathcal{C} = (C_{1}, \dots, C_{c})$ of $[0,1]^d$
such that for some $C_{i}$ the distribution $\mathcal{L}(\mathbf{U}
\mid
\mathbf{U} \in C_{i})$ is associated. Let $G$ be an increasing set and
let $\mathcal{B} =(C_1,\ldots,C_{i-1}, C_{i} \cap G, C_{i} \cap G^{c},
C_{i+1},\ldots,C_c)$.  If $f\dvtx[0,1]^d \rightarrow [0,1]$ is
non-decreasing, then $W_{\IE}^{\mathcal{B}} \le_{\cx}
W_{\IE}^{\mathcal{C}}$.
\end{lemma}

\begin{pf}
With $\mathcal{L}(\mathbf{V}_{1}) = \mathcal{L}(\mathbf{U}\mid\mathbf{U}
\in  C_{i} \cap G^c)$ and $\mathcal{L}(\mathbf{V}_{2}) =
\mathcal{L}(\mathbf{U}\mid\mathbf{U} \in  C_{i} \cap G)$,
Lemma~\ref{le:assoc} yields $\mathbf{V}_{1} \le_{\st} \mathbf{V}_{2}$.
The monotonicity of $f$ implies $f(\mathbf{V}_{1}) \le_{\st}
f(\mathbf{V}_{2})$, and Lemma~\ref{le:partmon1} now proves the claim,
applying arguments as in the proof of
Theorem~\ref{th:increasinginterror}.
\end{pf}

The following result can be found in \cite{KarRinIdsdJMVA1980}.

\begin{lemma}\label{le:MTP}
If an MTP$_2$ vector $\mathbf{U}$ takes values in a lattice of which
$C$ is a sublattice, then $\mathcal{L}(\mathbf{U} \mid \mathbf{U} \in C )$
is MTP$_2$ and hence associated.
\end{lemma}

The following corollary is obvious, and only requires the fact that the
intersection of sublattices is a lattice.

\begin{corollary}\label{cor:MTP}
If an MTP$_2$ vector $\mathbf{U}$ takes values in some lattice, and
$C$,  $G$ and $G^c$, are all sublattices, then both
$\mathcal{L}(\mathbf{U} \mid \mathbf{U} \in C \cap G)$  and
$\mathcal{L}(\mathbf{U} \mid \mathbf{U} \in C \cap G^c)$ are MTP$_2$, and
hence also associated.
\end{corollary}

\begin{pf*}{Proof of Theorem~\ref{th:multconvexintgen}}
We first prove by induction that $\mathcal{L}(\mathbf{U}\mid\mathbf{U} \in
C_{i,k})$ are MTP$_2$ for all $C_{i,k} \in \mathcal{C}_k$ and
$k=1,\ldots,m$. For $k=1$ this follows from Lemma \ref{le:MTP} and the
assumptions that $\mathbf{U}$ is MTP$_2$ and that $C_{i}=C_{i,1}$ are
sublattices of $[0,1]^d$. Assuming the statement true for $1 \le k<m$,
to verify that it is true for $k+1$ we need only show that
$\mathcal{L}(\mathbf{U}\mid\mathbf{U} \in C_{i_k,k}\cap G_k)$ and
$\mathcal{L}(\mathbf{U}\mid\mathbf{U} \in C_{i_k,k}\cap G_k^c)$ are
MTP$_2$, which follows from Lemma \ref{le:MTP}, thus completing the
induction.

Hence, again using Lemma \ref{le:MTP},
$\mathcal{L}(\mathbf{U}\mid\mathbf{U} \in C_{i_k,k})$ is associated. Since
$G_k$ is increasing, Lemma~\ref{le:multconvexint} now yields
$W_{\IE}^{\mathcal{C}_{k+1}} \le_{\cx} W_{\IE}^{\mathcal{C}_k}
 \mbox{ for all }k=1, \ldots, m-1$, and, therefore, the theorem.
\end{pf*}

\begin{figure}[b]

\includegraphics{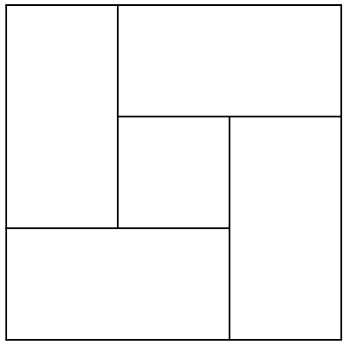}

\caption{Non-attainable tiling.}\label{fi:tiling}
\end{figure}

A sequence of partitions as in Theorem~\ref{th:multconvexintgen} can be
generated as follows:  start with the whole space $[0,1]^d$, then
split it into boxes by repeatedly  subdividing one element of the
partition by an intersection with some  $G$ and $G^c$. In $[0,1]^2$,
the resulting partition forms a tiling of the square by rectangles.
Note that from the first step, a sequence of partitions created using
$G$ as above has at least one line that crosses the whole square from
side to side. Therefore the tiling of Figure~\ref{fi:tiling} is not
attainable by such a sequence.

Finally, recall that the hypothesis of MTP$_{2}$  includes as a
particular case the uniform distribution on $[0,1]^{d}$, so
Theorem~\ref{th:multconvexintgen} applies to the estimation of the
integral $\int f(\mathbf{u})\,\diff{\mathbf{u}}$ on $[0,1]^{d}$, or any
lattice.

\begin{appendix}\label{se:appendix}

\section*{Appendix}
\setcounter{equation}{0}
\renewcommand{\thelemmasss}{A.1}
\begin{lemmasss}\label{le:supremumXi}
Given a partition $\mathcal{B}^{*}$ of $N$, consider a collection of
independent random variables $\{\xi^{B^{*}}_{j}\}$,  $B^{*} \in
\mathcal{B}^{*}$, $j\in {B}^{*}$, with those indexed by  the same
element $B^{*}$ of the partition being identically distributed.

For $\mathcal{C}^{*} \le_{\refin} \mathcal{B}^{*}$, let
$\{\xi^{C^{*}}_{j}\}$ with $C^{*} \in \mathcal{C}^{*}$ and $j\in C^{*}$
be a collection of independent random variables with the mixture
distribution
\begin{equation}\label{eq:LxiBC}
\mathcal{L}(\xi^{C^{*}}_{j}) = \sum_{B^{*} \subset C^{*}}
\frac{|B^{*}|}{|C^{*}|} \mathcal{L}(\xi^{B^{*}}_{j}).
\vspace{-3pt}
\end{equation}

Then
\begin{equation}\label{eq:maxxi}
\vspace{-3pt}
\max_{C^{*} \in \mathcal{C}^{*}}\max_{j \in C^{*}} \xi^{C^{*}}_{j}
\le_{\st} \max_{B^{*} \in \mathcal{B}^{*}}\max_{j \in B^{*}}
\xi^{B^{*}}_{j}.
\end{equation}
\end{lemmasss}

\begin{pf}
Let $ p^{B^{*}} = \mathbb{P}(\xi^{B^{*}}_{1} \le t)$ for $B^{*} \in
\mathcal{B}^{*}$ and $p^{C^{*}} = \mathbb{P}(\xi^{C^{*}}_{1} \le t)$
for $C^{*} \in \mathcal{C}^{*}$.

We claim that
\[
\vspace{-3pt}
(\underbrace{p^{C_{1}^{*}},\ldots,p^{C_{1}^{*}}}_{|C_{1}^{*}|}, \dots,
\underbrace{p^{C_{c}^{*}},\ldots,p^{C_{c}^{*}}}_{|C_{c}^{*}|}) \prec
(\underbrace{p^{B_{1}^{*}},\ldots,p^{B_{1}^{*}}}_{|B_{1}^{*}|}, \dots,
\underbrace{p^{B_{b}^{*}},\ldots,p^{B_{b}^{*}}}_{|B_{b}^{*}|}).
\]
To see this, observe that \eqref{eq:LxiBC} implies that the vector on
the left-hand side above is obtained from the one on the right by
multiplying it by the $n \times n$ doubly stochastic matrix
$\mathbf{D}$, which is block diagonal where the $i$th block is the
$|C_{i}^{*}| \times |C_{i}^{*}|$ matrix with all entries equal to
$1/|C_{i}^{*}|$.

Hence, by the Schur concavity of the function $(\theta_{1}, \dots,
\theta_{n}) \mapsto \prod_{i=1}^n \theta_i$, we have
\[
\mathbb{P}\Bigl(\max_{C^{*} \in \mathcal{C}^{*}}\max_{j \in C^{*}}
\xi^{C^{*}}_{j} \le t \Bigr) = \prod_{C^{*} \in \mathcal{C}^{*}}
(p^{C^{*}})^{|C^{*}|} \ge \prod_{B^{*} \in \mathcal{B}^{*}}
(p^{B^{*}})^{|B^{*}|} = \mathbb{P}\Bigl(\max_{B^{*} \in
\mathcal{B}^{*}}\max_{j \in B^{*}} \xi^{B^{*}}_{j} \le t\Bigr),
\]
which is equivalent to \eqref{eq:maxxi}.
\end{pf}

\begin{pf*}{Proof of Theorem~\ref{th:maxgeneral}}
Let $\mathcal{B}^{*}$ and $\mathcal{C}^{*}$ be partitions associated
with $\mathcal{B}$ and $\mathcal{C}$, respectively, satisfying
$\mathcal{C}^{*} \le_{\refin} \mathcal{B}^{*}$, and let
$\{\xi^{B^{*}}_{j},B^{*} \in \mathcal{B}^{*}, j \in B^{*} \}$ and
$\{\xi^{C^{*}}_{j},C^{*} \in \mathcal{C}^{*}, j \in C^{*} \}$ be
collections of independent random variables with distributions\vspace{-5pt}
\begin{eqnarray*}
\mathbb{P}(\xi^{B^{*}}_{j} \le t)
&=&
\mathbb{P}\bigl(f(U)\le t\mid U \in B\bigr),
\\
\mathbb{P}(\xi^{C^{*}}_{j} \le t)
&=&
\mathbb{P}\bigl(f(U)\le t\mid U \in C\bigr).
\end{eqnarray*}
\vspace{-5pt}
Then  \eqref{eq:LxiBC} holds (law of total probability), and the result
follows by Lemma~\ref{le:supremumXi}.
\end{pf*}

\begin{pf*}{Proof of Theorem~\ref{th:interror}}
In what follows we consider conditional expectation with
respect to a partition. Though the notion is standard, specifically, by
$\mathbb{E}[f(U) + \varepsilon \vert \mathcal{B}]$, we mean the random
variable that takes values $\overline{f}_{B} := \mathbb{E}[f(U)\mid U \in
B]$ with probability $|B^{*}|/n$. Then
\begin{eqnarray*}
\Var [f(U) + \varepsilon \vert \mathcal{B}]
&=&
\mathbb{E}\bigl[\{f(U) +\varepsilon - \mathbb{E}[f(U)+ \varepsilon \vert\mathcal{B}]\}^2 \vert\mathcal{B}\bigr]
\\
&=&
\mathbb{E}\bigl[\{f(U) + \varepsilon - \mathbb{E}[f(U)\vert\mathcal{B}]\}^2\vert\mathcal{B}\bigr]
\end{eqnarray*}
is a random variable taking  values $\mathbb{E}[(f(U) + \varepsilon
-\overline{f}_{B})^2 \mid U \in B]$  with probability $|B^{*}|/n$, and
\begin{eqnarray*}
\mathbb{E}\bigl[\Var [f(U) + \varepsilon \vert \mathcal{B}]\bigr]
&=&
\sum_{B\in \mathcal{B}}\frac{|B^{*}|}{n} \mathbb{E}\bigl[\bigl(f(U) + \varepsilon -\overline{f}_{B}\bigr)^2\mid U \in B\bigr] \\
&=&
\frac{1}{n}\sum_{B\in \mathcal{B}} |B^{*}| \mathbb{E}\bigl[\bigl(f(V_1^B) + \varepsilon -\overline{f}_{B}\bigr)^2\bigr] \\
&=&
\frac{1}{n} \Var\biggl[\sum_{B\in \mathcal{B}} \sum_{j \in B_{i}^{*}} f(V_{j}^{B}) + \varepsilon^{B}_{j}\biggr]\\
&=&
n \Var[W_{\IE}^{\mathcal{B}}].
\end{eqnarray*}

If $ \mathcal{C} \le_{\refin}\mathcal{B}$, then for any random variable
$Y$, say, $\Var [\mathbb{E}[Y | \mathcal{B}]] \ge \Var [\mathbb{E}[Y |
\mathcal{C}]]$ by Jensen's inequality, and now the usual variance
decomposition of $Y$ (see, e.g., \cite{RosdsfWorld2006}, Theorem~13.3.1)
implies $ \mathbb{E}[\Var [Y | \mathcal{B}]] \le \mathbb{E}[\Var [Y |
\mathcal{C}]]$. Therefore
\[
\mathbb{E}\bigl[\Var [f(U) + \varepsilon \vert \mathcal{B}]\bigr] \le \mathbb{E}\bigl[\Var
[f(U) + \varepsilon \vert \mathcal{C}]\bigr],
\]
and hence
\[
\Var[W_{\IE}^{\mathcal{B}}] = \frac{1}{n}\mathbb{E}\bigl[\Var [f(U) +
\varepsilon \vert \mathcal{B}]\bigr] \le \frac{1}{n}\mathbb{E}\bigl[\Var [f(U) +
\varepsilon \vert \mathcal{C}]\bigr] = \Var[W_{\IE}^{\mathcal{C}}].
\]
\upqed
\end{pf*}
\end{appendix}

\section*{Acknowledgements}
We thank Abram Kagan for sparking our curiosity in the topic with a
simple version of Theorem~\ref{th:maxgeneral}, Erich Novak for an
important bibliographical reference, and Pierpaolo Brutti for his help
with R. We are indebted to the editor, an associate editor and three
referees for their accurate reading of the paper and their helpful
comments.
The work of Yosef Rinott is partially supported by the Israel Science
Foundation grant No. 473/04. The work of Marco Scarsini is partially
supported by  MIUR-COFIN.

\printhistory


\begin{thebibliography}{00}

\bibitem{BaiDurdsdJSPI1997}
Bai, S.K. and Durairajan, T.M. (1997). Optimal equivariant estimator
with respect to convex loss function. \textit{J. Statist.
Plann. Inference} \textbf{64} 283--295.
\MR{1621618}


\bibitem{BerdsdAS1976}
Berger, J.O. (1976). Admissibility results for generalized {B}ayes
estimators of coordinates of a location vector.
\textit{Ann. Statist.} \textbf{4} 334--356.
\MR{0400486}


\bibitem{BladsdP2BSMSP1951}
Blackwell, D. (1951). Comparison of experiments. In
\textit{Proceedings of the {S}econd {B}erkeley {S}ymposium on
{M}athematical {S}tatistics and {P}robability, 1950} 93--102. Berkeley
and Los Angeles, CA: California Univ. Press.
\MR{0046002}


\bibitem{BladsdAMS1953}
Blackwell, D. (1953). Equivalent comparisons of experiments.
\textit{Ann. Math. Statist.} \textbf{24} \mbox{265--272}.
\MR{0056251}


\bibitem{EbedsdSD1984}
Eberl Jr., W. (1984). On unbiased estimation with convex loss
functions. \textit{Statist. Decisions} \textbf{1984} 177--192.
\MR{0785208}


\bibitem{ErmZhiKondsdDANSSSR1988}
Ermakov, S.M., Zhiglyavski{\u\i}, A.A. and Kondratovich, M.V. (1988).
Reduction of a problem of random estimation of an extremum of a
function. \textit{Dokl. Akad. Nauk SSSR} \textbf{302} \mbox{796--798}.
\MR{0983943}


\bibitem{GladsdSpringer2004}
Glasserman, P. (2004). \textit{Monte {C}arlo Methods in Financial
Engineering}. New York: Springer.
\MR{1999614}


\bibitem{GolRinScadsdmimeo2010}
Goldstein, L., Rinott, Y. and Scarsini, M. (2010). Stochastic
comparisons of stratified sampling techniques for some Monte Carlo
estimators. Technical report. Available at
\href{http://arxiv.org/abs/1005.5414v1}{arXiv:1005.5414v1} [math.ST].


\bibitem{KarNovdsdPJM1963}
Karlin, S. and Novikoff, A. (1963). Generalized convex inequalities.
\textit{Pacific J. Math.} \textbf{13} 1251--1279.
\MR{0156927}


\bibitem{KarRinIdsdJMVA1980}
Karlin, S. and Rinott, Y. (1980). Classes of orderings of measures and
related correlation inequalities. {I}. {M}ultivariate totally positive
distributions. \textit{J. Multivariate Anal.} \textbf{10} \mbox{467--498}.
\MR{0599685}


\bibitem{KonZhidsdLNS1998}
Kondratovich, M. and Zhigljavsky, A. (1998). Comparison of independent
and stratified sampling schemes in problems of global optimization. In
\textit{Monte {C}arlo and Quasi-{M}onte {C}arlo Methods 1996
({S}alzburg)} 292--299. New York: Springer.
\MR{1644527}


\bibitem{KozdsdJMVA1977}
Kozek, A. (1977). Efficiency and {C}ram\'er--{R}ao type inequalities
for convex loss functions. \textit{J. Multivariate Anal.} \textbf{7}
89--106.
\MR{0431482}


\bibitem{LaydsdJRSSB1972}
Laycock, P.J. (1972). Convex loss applied to design in regression
problems. \textit{J. Roy. Statist. Soc. Ser. B} \textbf{34} 148--170,
170--186.
\MR{0350935}


\bibitem{LaySildsdB1968}
Laycock, P.J. and Silvey, S.D. (1968). Optimal designs in regression
problems with a general convex loss function. \textit{Biometrika}
\textbf{55} 53--66.
\MR{0225446}


\bibitem{LinMoudsdAISM1982}
Lin, P.E. and Mousa, A. (1982). Proper {B}ayes minimax estimators for a
multivariate normal mean with unknown common variance under a convex
loss function. \textit{Ann. Inst. Statist. Math.} \textbf{34} 441--456.
\MR{0695065}


\bibitem{MarOlkdsdACADEMIC1979}
Marshall, A.W. and Olkin, I. (1979). \textit{Inequalities: Theory of
Majorization and Its Applications}. New York: Academic Press.
\MR{0552278}


\bibitem{MueStodsdWiley2002}
M{\"u}ller, A. and Stoyan, D. (2002). \textit{Comparison Methods for
Stochastic Models and Risks}. Chichester: Wiley.
\MR{1889865}


\bibitem{NovdsdLNMSpringer1988}
Novak, E. (1988). \textit{Deterministic and Stochastic Error Bounds in
Numerical Analysis}. Berlin: Springer.
\MR{0971255}


\bibitem{PapdsdJC1993}
Papageorgiou, A. (1993). Integration of monotone functions of several
variables. \textit{J. Complexity} \textbf{9} 252--268.
\MR{1226312}


\bibitem{PetKoureAISM2001}
Petropoulos, C. and Kourouklis, S. (2001). Estimation of an exponential
quantile under a general loss and an alternative estimator under
quadratic loss. \textit{Ann. Inst. Statist. Math.} \textbf{53}
746--759.
\MR{1880809}


\bibitem{RosdsfWorld2006}
Rosenthal, J.S. (2006). \textit{A First Look at Rigorous Probability
Theory}, 2nd ed. Hackensack, NJ: World Scientific Publishing.
\MR{1767078}


\bibitem{RosSchfrdeOR1984}
Ross, S.M. and Schechner, Z. (1984). Some reliability applications of
the variability ordering. \textit{Oper. Res.} \textbf{32} 679--687.
\MR{0756013}


\bibitem{ShafdJMVA1982}
Shaked, M. (1982). A general theory of some positive dependence
notions. \textit{J. Multivariate Anal.} \textbf{12} 199--218.
\MR{0661559}


\bibitem{ShaShadsfSPRI2007}
Shaked, M. and Shanthikumar, J.G. (2007). \textit{Stochastic Orders}.
New York: Springer.
\MR{2265633}


\bibitem{ZhiZilfdSPRINGER2008}
Zhigljavsky, A. and {\v{Z}}ilinskas, A. (2008). \textit{Stochastic
Global Optimization}. New York: Springer.
\MR{2361744}


\bibitem{ZhiChedsMCM1996}
Zhigljavsky, A.A. and Chekmasov, M.V. (1996). Comparison of
independent, stratified and random covering sample schemes in
optimization problems. \textit{Math. Comput. Modelling} \textbf{23}
97--110.
\MR{1398005}


\end{thebibliography}
\end{document}